\documentclass[a4paper, 11pt]{amsart}   
\usepackage{mathptmx, amssymb,amscd,latexsym, eulervm,calligra,mathrsfs,amsmath,amsthm,amsfonts,amsrefs,tikz}
\usepackage[colorlinks=true,linkcolor=blue,urlcolor=blue,citecolor=blue]{hyperref}
\usepackage{xurl}
\hypersetup{breaklinks=true}
\usepackage{setspace}
\usepackage{tabularx}
\usepackage{comment}
\usepackage{paralist}
\usepackage{aliascnt}
\usepackage{blkarray}
\usepackage{booktabs}
\usepackage{enumitem}
\setlist{leftmargin=5.5mm}
\usepackage[inner=2.5cm,outer=2.5cm, bottom=3.4cm]{geometry}
\usepackage[capitalize,nameinlink]{cleveref}
\usepackage{tikz, tikz-cd}
\usepackage{calligra,mathrsfs}
\usetikzlibrary{matrix}
\usetikzlibrary{arrows,calc}
\usepackage{url}
\usepackage{amssymb,latexsym,amsmath,amsfonts,amsthm,amscd,graphicx,url,color,hyperref,stmaryrd}

\theoremstyle{plain}

\newtheorem{thm}{Theorem}

\theoremstyle{definition}

\newtheorem{problem}{Problem}
\newtheorem{solved}[problem]{(Solved) Problem}

\newcommand{\scR}{\mathscr{R}}
\newcommand{\scT}{\mathscr{T}}
\newcommand{\scC}{\mathscr{C}}
\newcommand{\scD}{\mathscr{D}}

\newcommand{\G}{\Gamma}
\newcommand{\edim}{\mathrm{edim}}
\newcommand{\mult}{\mathrm{mult}}
\newcommand{\type}{\mathrm{t}}
\newcommand{\reg}{\mathrm{reg}}
\newcommand{\width}{\mathrm{width}}
\newcommand{\Frob}{\mathrm{Frob}}
\newcommand{\tG}{{\tilde{\Gamma}}}
\newcommand{\mm}{\mathfrak{m}}
\newcommand{\NN}{\mathbb{N}}
\newcommand{\ZZ}{\mathbb{Z}}

\newcommand{\gr}{\mathrm{gr}}

\newcommand{\kk}{\Bbbk}

\makeatletter
\@namedef{subjclassname@2020}{
  \textup{2020} Mathematics Subject Classification}
\makeatother

\begin{document}

\author[A.\,Moscariello, A.\,Sammartano]{Alessio~Moscariello and Alessio~Sammartano}
\address{Alessio Moscariello: Dipartimento di Matematica e Informatica\\Università degli Studi di Catania\\Catania\\Italy}
\email{alessio.moscariello@unict.it}
\address{Alessio Sammartano: Dipartimento di Matematica \\ Politecnico di Milano \\ Milan \\ Italy}
\email{alessio.sammartano@polimi.it}

\title{Open problems on relations of numerical semigroups}


\begin{abstract}
In this paper, 
we collect some open problems about  minimal presentations of numerical semigroups and, more generally, 
about defining ideals and free resolutions of their semigroup rings and associated graded rings. 
We emphasize both long-standing problems and more recent questions and developments.

Added 2026:  slightly updated version of the problem list,  marking the problems that were solved.
\end{abstract}

\maketitle


\section*{Introduction}

A classical goal of algebra is to study an object by means of generators and relations among the generators.
In this paper, we consider this problem for some of the most elementary objects in mathematics: 
  \emph{numerical semigroups}, i.e. (cofinite) additive submonoids $\G \subseteq \mathbb{N}$.
Despite the apparent simplicity, 
this topic is surprisingly rich, and ripe with  interesting connections with various areas of mathematics,
such as  combinatorics, number theory, algebraic geometry, and coding theory.
Our point of view is mainly that of commutative algebra, and some  questions in this paper  are inspired by this perspective.

The purpose of this paper is to 
highlight the richness of this area, by
advertising some of the most interesting open questions about relations of numerical semigroups.\footnote{Some previous surveys on related topics are
\cite{Barucci,BGT,GS20,Stamate}.}
We interpret the term ``relations'' in a broad sense.
The primary object of interest is minimal presentations of a numerical semigroup $\G$, 
that is, minimal collections of relations among the generators of $\G$.
However,
this is equivalent to studying the generators of the toric ideal of the local semigroup ring $\kk\llbracket \G\rrbracket$,
or of the graded semigroup ring $\kk[\G]$.
This simple observation points to several natural and interesting generalizations, such as
the investigation of the higher syzygies of the semigroup rings, 
 of the defining ideal of the associated graded ring
$\gr(\kk\llbracket \G\rrbracket)$, 
or of more general 1-dimensional local rings.

A secondary aim of this paper is to emphasize the advantages brought by 
casting  problems in the framework of commutative algebra, in particular, 
in the framework of free resolutions.
One is the uniform treatment of several different invariants of a numerical semigroup $\G$:
for example,
the number of minimal relations and the type of $\G$ are respectively the first and last Betti numbers of the semigroup rings, while the Frobenius number of $\G$ can be read from the degrees of the syzygies of $\kk[\G]$.
Another advantage is that conditions like being a symmetric or complete intersection numerical semigroup assume a clear interpretation in terms of  free resolutions.
Finally, 
a main advantage of this approach  is the access it gives to a variety of tools developed in commutative algebra  for the analysis of Hilbert functions, Betti numbers, Hilbert schemes of points, etc.

\tableofcontents

\section{Generalities}\label{SectionSetup}

In this section, we fix the basic notation for this paper and introduce the main concepts related to numerical semigroups and the various commutative rings associated to them. 
We refer to 
\cite{RosalesGarciaSanchezBook}
for background on numerical semigroups, and to \cite{Eisenbud,HerzogHibi} for 
commutative algebra notions.

Let $\NN$ denote the set of non-negative integers.
A numerical semigroup  $\G\subseteq \NN$ is a cofinite additive submonoid of $\NN$.
The cofiniteness condition is often unnecessary, since every non-zero additive submonoid of $\NN$ is isomorphic to a numerical semigroup.
A numerical semigroup has a unique minimal set of  generators $ g_1 < g_2 < \cdots < g_e \in \G$, 
in this case we write $\G= \langle  g_1, \ldots, g_e\rangle$.
The number $e$ of  minimal generators is called the embedding dimension of $\G$, 
denoted by $\edim(\G)$.
The smallest positive integer $g_1$ belonging to $\G$ is called the multiplicity of $\G$, denoted by $\mult(\G)$.
An integer $p \in \ZZ\setminus \G$ is called a pseudo-Frobenius number of $\G$ if $p+\gamma \in \G$ for all $\gamma \in \G\setminus\{0\}$.
The set of pseudo-Frobenius numbers of $\G$ is finite, and its cardinality is called the type of $\G$, denoted by $\type(\G)$.
We have the inequalities $\edim(\G) \leq \mult(\G)$ and $\type(\G) \leq \mult(\G) -1$.
The largest pseudo-Frobenius number of $\G$ is called the Frobenius number of $\G$, denoted by $\Frob(\G)$.

There is a semigroup homomorphism  $\varphi: \NN^{e} \to \NN$ defined by $\varphi( a_1 , \ldots, a_e) = \sum_{i=1}^e a_i g_i$.
The kernel congruence of $\varphi$ is  $\ker(\varphi) = \{(\mathbf{a},\mathbf{b}) \in \NN^{e} \times \NN^{e} \, \mid \, \varphi(\mathbf{a}) = \varphi(\mathbf{b})\}$,
and a minimal set of generators of the congruence $\ker(\varphi)$ is called a
 minimal presentation of $\G$.
 The cardinality of a minimal presentation of $\G$ depends only on $\G$ (cf. \cite[Corollary 8.13]{RosalesGarciaSanchezBook}),  we denote it by  $\rho(\G)$ and 
 refer to it informally as  the \emph{number of relations} of $\G$.

Let $\kk$ be a field. 
We now describe various ways to associate a commutative  $\kk$-algebra to a numerical semigroup $\G$.
Each of them has its advantages and disadvantages, and can be used in different ways to apply commutative algebra methods to  obtain information on the semigroup.

The {\bf graded semigroup ring} is $A_\G = \kk[\G] = \kk[\tau^\gamma \, \mid \, \gamma \in \G] \subseteq \kk[\tau]$, where $\tau$ is a variable.
It is a 1-dimensional positively graded Noetherian $\kk$-domain, where $\deg(\tau^\gamma)= \gamma \in \NN$.
In particular, it is Cohen-Macaulay.
Notice that $A_\G$ is not a standard graded $\kk$-algebra unless $\G = \NN$, since its algebra generators have degrees greater than 1.

The {\bf local semigroup ring} is $R_\G  = \kk\llbracket\G\rrbracket = \kk\llbracket \tau^\gamma \, \mid \, \gamma \in \G\rrbracket \subseteq \kk\llbracket \tau\rrbracket$, that is, the completion of $A_\G$ with respect to its homogeneous maximal ideal. 
It is a 1-dimensional Noetherian local domain, and it is Cohen-Macaulay.
Its multiplicity as a local ring coincides with $\mult(\G)$.
For most purposes, one could consider instead the localization of $A_\G$ at its homogeneous maximal ideal. 

The {\bf associated graded ring} of $R_\G$ is $G_\G = \gr(R_\G) = \bigoplus_{i \geq 0} {\mm^{i}}/{\mm^{i+1}}$,
where $\mm$ denotes the  maximal ideal of $R_\G$.
It is a 1-dimensional standard graded Noetherian $\kk$-algebra.
However, it is not a domain unless $\G= \NN$, and often it is not even Cohen-Macaulay.
In all three cases $A_\G, R_\G, G_\G$, 
the embedding dimension of the ring, i.e. the number of generators of the maximal ideal, coincides with $\edim(\G)$.

Let $R$ be any Noetherian local (respectively, positively graded) $\kk$-algebra. 
A minimal regular presentation of $R$ is a surjective local (respectively, graded) $\kk$-algebra homomorphism $S \to R$ such that $S$ is a regular local ring (respectively, a graded polynomial ring over $\kk$) and such that $\edim(R) = \edim(S)$. 
The  Betti numbers of $R$ over $S$ are
the integers 
$
b_i(R) = \dim_\Bbbk \mathrm{Tor}_i^S(R,\Bbbk).
$
They are the ranks of the free modules in a minimal free resolution of $R$ as $S$-module.
The notation  is justified, since they are independent of the choice of a minimal regular presentation.
We have $b_0(R) = 1$, whereas $b_1(R)$ is the number of generators of the defining ideal of $R$,
that is, the kernel of the presentation $S \to R$.
Moreover, 
\begin{equation}\label{EqAlternating}
\sum_{i\geq 0} (-1)^ib_i(R) =0.
\end{equation}

If $R$ is graded, then we can also consider the graded Betti numbers of $R$ over $S$, defined as 
$
b_{i,j}(R) = \dim_\Bbbk [\mathrm{Tor}_i^S(R,\Bbbk)]_j.
$
Here,
the notation $[V]_j$ indicates the graded component of degree $j$ of  a graded vector space $V$.
The integer $b_{i,j}(R)$ counts the number of $i$th syzygies of degree $j$, that is, 
the number of elements of degree $j$ in a basis of the $i$th  free module in a minimal graded free resolution of $R$ over $S$.
Clearly, we have $\sum_{j} b_{i,j}(R) = b_i(R)$ and, thus, the Betti numbers $b_i(R)$ are sometimes called the \emph{total} Betti numbers, to avoid confusion with the graded ones.

Let $\G$ be a numerical semigroup.
We have $b_i(A_\G) = b_i(R_\G)\leq b_i(G_\G)$ for all $i$.
Moreover, we have $b_0(A_\G)=1, b_1(A_\G) = \rho(\G), b_{e-1}(A_\G)=\type(\G)$, and $b_i(A_\G) \ne 0$ if and only if 
$0 \leq i \leq e-1$.
On the other hand, we have $b_i(G_\G) = 0$  if $ i > e$, 
and $b_{e}(G_\G) = 0$  if and only if $G_\G$ is Cohen-Macaulay.

The  Betti numbers (graded or total) of $A_\G$ can be computed in terms of homology groups of certain simplicial complexes associated to $\G$,
cf. \cite{BrunsHerzog,CampilloMarijuan}.
It follows from \cite[Theorem 2.1]{BrunsHerzog} that, in general, they depend on the characteristic of the ground field.
(Observe that
the first and last Betti numbers do not depend on the characteristic, as they  can be read off from combinatorial data of $\G$.)
It is unknown whether a simplicial description exists for the associated graded ring.

\begin{problem}\label{ProblemSimplicialAssociatedGraded}
Let $\G$ be a numerical semigroup. 
Is it possible to compute 
the Betti numbers $b_i(G_\G)$ or $b_{i,j}(G_\G)$ via homology of appropriate simplicial complexes?
\end{problem}

\section{Embedding dimension and multiplicity}\label{SectionEdimMult}

When studying  Betti numbers of semigroup rings, it is natural to classify numerical semigroups by  their embedding dimension.
The only semigroup with $\edim(\G) =1$ is $\G= \NN$.
When $\edim(\G) =2$,  $R_\G$ is the coordinate ring of an affine plane curve, and $b_1(\G) =1$.
The case $\edim(\G)=3$ is completely described in the following theorem.

\begin{thm}[\protect{\cite{Herzog}}]\label{ThmHerzog}
Let $\G$ be a numerical semigroup with $\edim(\G)= 3$.
If $\G$ is symmetric, then $b_1(R_\G) =2 $ and $b_2(R_\G) =1$.
If $\G$ is not symmetric, then $b_1(R_\G) =3 $ and $b_2(R_\G) =2$.
\end{thm}

See Section \ref{SectionSymmetric} for the definition of symmetric semigroup.
The situation becomes more complicated when $\edim(\G)\geq 4$.
Bresinsky \cite{Bresinsky75b} discovered a family of semigroups of embedding dimension 4 with arbitrarily large number of minimal relations,
while 
Backelin \cite[p.75]{FGH} discovered a family of semigroups of embedding dimension 4 with arbitrarily large type.
Another  such family appeared in \cite{Arslan}.
By applying the gluing construction \cite{Gluing}, 
one can easily produce families with the same behavior and $\edim(\G) > 4$.

It is, therefore, necessary to refine the analysis by taking into account other invariants besides the embedding dimension. 
The most natural choice is the multiplicity,
since Betti numbers of Cohen-Macaulay rings of bounded multiplicity are bounded.
Specifically, if $R$ is a Cohen-Macaulay ring with multiplicity $m$, we have
$b_i(R) \leq i{m \choose i+1}$ for all $1 \leq i \leq m-1$,
 see for example \cite{Valla}.
 The extremal case is well-understood:
 it is not hard to see that the following conditions are equivalent.
\begin{itemize}
\item $\G$ has maximal embedding dimension $\edim(\G) = \mult(\G)$;
\item $\G$ has maximal type  $\type(\G) = \mult(\G)-1$;
\item $R_\G$ has maximal Betti number $b_i(R_\G) = i{\mult(\G)\choose i+1}$ for some, equivalently, every, $1 \leq i\leq \mult(\G)-1$;
\item $G_\G$ has maximal Betti number $b_i(G_\G) = i{\mult(\G)\choose i+1}$ for some, equivalently, every, $1 \leq i\leq \mult(\G)-1$.
\end{itemize}

The central problem in this section is the following.

\begin{problem}[\protect{\cite[Question 7.2]{EHOOPK}}]\label{ProblemRem}
Let $e, m \in \NN$ with $ 3 \leq e < m$.
Determine 
$$
\mathscr{R}(e,m) = \sup \big \{ \rho(\G) \, \mid \,  \edim(\G) = e+1,\, \mult(\G) = m\big\}.
$$
\end{problem}

When  $m-e$ is small, some values of  $\mathscr{R}(e,m)$ were determined in the following theorem.

\begin{thm}[\protect{\cite{RosalesGarciaSanchez98}}]\label{TheoremRem12}
If $m -e \leq 3$, then 
$ \scR(e,m) ={e+1 \choose 2}$.
\end{thm}

This result was later extended in \cite{MoscarielloSammartano}.

\begin{thm}[\protect{\cite[Theorem 3.3]{MoscarielloSammartano}}]\label{TheoremRem345}
We have $\mathscr{R}(e,e+4)=\mathscr{R}(e,e+5)= {e+1 \choose 2}+1$ and 
$\mathscr{R}(e,e+6)={e+1 \choose 2}+2$.
\end{thm}

As observed  in \cite[Section 1.3]{RosalesGarciaSanchez98},
Problem \ref{ProblemRem} becomes more complicated as $m-e$ increases.
In particular, for each fixed value of  $\delta = m-e\geq 6$, unlike Theorems \ref{TheoremRem12} and \ref{TheoremRem345},
it is not possible to obtain a polynomial formula for $\mathscr{R}(e,e+\delta)$ that is valid for all  $e$.
For example, we have the following result.

\begin{thm}[\protect{\cite[Theorem 3.2, Example 4.3, Proposition 4.6]{MoscarielloSammartano}}]\label{TheoremRem6}
We have $\mathscr{R}(e,e+7) = {e+1 \choose 2}+4$ for $e\geq 9$, but  $\mathscr{R}(e,e+7) < {e+1 \choose 2}+4$ for $e = 3,4,5$.
\end{thm}

If we consider instead the smallest value of the embedding dimension for which the problem is open,
the following conjecture proposes an answer to Problem \ref{ProblemRem}.

\begin{problem}[\protect{\cite[Conjecture 7.1]{EHOOPK}}]\label{ProblemR4m}
We have $
\mathscr{R}(3,m) = \lfloor 2 + 2\sqrt{m}\rfloor.
$
\end{problem}

Now, we discuss a general upper bound for $\mathscr{R}(e,m)$ obtained in \cite{ERV}.
For each  $n,d \in \NN_{>0}$,
there exist unique integers  $n_d > n_{d-1} > \cdots > n_j \geq j \geq 1$
such that  $n = {n_d \choose d} + {n_{d-1}\choose d-1} + \cdots + {n_j \choose j}$. 
Let $n^{\langle d\rangle}= \binom{n_d+1}{d+1}+\binom{n_{d-1}+1}{d}+\cdots+\binom{n_j+1}{j+1}.$
Let $e,m \in \NN$ be such that $3 \leq e < m$,  let $r$ be the integer such that ${e+r-1 \choose r-1} \leq m < {e+r \choose r}$,
and let $s = m-{e+r-1 \choose r-1}$.
 Define $\scC(e,m) = {e+r-1 \choose r} +s^{\langle r\rangle}-s.$

\begin{thm}[\protect{\cite[Theorem 5.1]{ERV}}]\label{TheoremERV}
We have $\scR(e,m) \leq \scC(e,m)$ for all $e,m$.
\end{thm}

In fact, \cite[Theorem 5.1]{ERV} applies to all Cohen-Macaulay local  rings, and Theorem \ref{TheoremERV} is the special case of $R_\G$.
Thus, a natural first step towards solving Problem \ref{ProblemRem} is the following:

\begin{problem}\label{ProblemRemCem}
Determine the values of $e,m$ for which $\scR(e,m) = \scC(e,m)$.
\end{problem}

A partial answer is provided by the following result.

\begin{thm}[\protect{\cite[Theorem 1.2]{MoscarielloSammartano}}]\label{TheoremRC}
Let $e,m \in \NN$ be such that $3 \leq e < m$.
\begin{enumerate}
\item For each fixed value of $\delta $, we have $\scR(e,e+\delta) = \scC(e,e+\delta)$ for all $e \gg 0$.
\item For each fixed value of $e$, we have $\scR(e,m) < \scC(e,m)$ for all $m \gg 0$.
\end{enumerate}
\end{thm}

One may also consider the analogue question for the type of numerical semigroups.

\begin{problem}\label{ProblemTem}
Determine 
$
\mathscr{T}(e,m) = \sup \big \{ \type(\G) \, \mid \,  \edim(\G) = e+1,\, \mult(\G) = m\big\}.
$
\end{problem}

In fact, the situation  for  $\mathscr{T}(e,m)$ is analogous to that of  $\scR(e,m)$.
Defining  the integers
$n_{\langle d\rangle}= \binom{n_d-1}{d}+\binom{n_{d-1}-1}{d-1}+\cdots+\binom{n_j-1}{j}$
and $\scD(e,m) = {e+r-2 \choose r-1} +s_{\langle r\rangle}$,
we have  $\scT(e,m) \leq \scD(e,m)$ by  \cite[Theorem 4]{EGV}. 
The following analogue of the first part of Theorem \ref{TheoremRC} holds:

\begin{thm}[\protect{\cite[Theorem 1.4]{MoscarielloSammartano}}]\label{TheoremTD}
Let $e,m \in \NN$ be such that $3 \leq e < m$.
 For each fixed value of $\delta$, we have $\scT(e,e+\delta) = \scD(e,e+\delta)$ for all $e \gg 0$.
\end{thm}

On the other hand, the analogue of the second part of Theorem \ref{TheoremRC} is not known:

\begin{problem}[\protect{\cite[Conjecture 5.3]{MoscarielloSammartano}}]\label{ProblemConjTem}
For each fixed value of $e \ge 3$, we have $\scT(e,m) < \scD(e,m)$ for all $m \gg 0$.
\end{problem}

Finally, it is possible to generalize the discussion to all the Betti numbers of $R_\G$, and obtain analogous results to 
Theorems \ref{TheoremRC} and \ref{TheoremTD},
cf. \cite{Valla,MoscarielloSammartano}.

\section{Width}\label{SectionWidth}

The  \emph{width} of a numerical semigroup $\G = \langle g_1 < g_2 < \cdots < g_e\rangle$ is  defined as  $\width(\G) = g_e-g_1$.
The theme of this section is finding bounds for the number of relations in terms of the width.

\begin{problem}\label{ProblemWidth}
Let $w \in \NN$. 
Determine $\mathscr{W}(w) := \sup\big\{\rho(\G) \,\mid \, \width(\G) = w\big\} $.
\end{problem}

It is not clear from the definition that $\mathscr{W}(w) < \infty $.
In fact, based on the discussions in Section \ref{SectionEdimMult},
one may perhaps expect  that $\mathscr{W}(w) = \infty $, 
since the multiplicity is allowed to be arbitrarily large.
Nevertheless, it follows from results of \cite{Vu} that $\mathscr{W}(w) < \infty $ for every $w \in \NN$.

When studying numerical semigroups with fixed width, a natural construction to consider is 
the interval completion of $\G$, defined as the semigroup 
$\tG = \langle g_1, g_1+1, g_1+2, \ldots, g_e\rangle$.
Observe that the given generating set is minimal if and only if $g_e < 2g_1$.
Moreover,  $\tG$ satisfies $\mult(\tG) = \mult(\G)$
and $\width(\tG) \leq \width(\G)$,
with $\width(\tG) = \width(\G)$
 if and only if  $g_e < 2g_1$.

In \cite{HerzogStamate}, Herzog and Stamate proposed the following conjectures.

\begin{problem}[\protect{\cite[Conjecture 2.4]{HerzogStamate}}]\label{ProblemHSgr1}
Let $\G$ be a numerical semigroup, then 
$b_1(G_\G) \leq b_1(G_\tG)$.
\end{problem}

\begin{problem}[\protect{\cite[Conjecture 2.1]{HerzogStamate}}]\label{ProblemHSgr2}
Let $\G$ be a numerical semigroup, then 
$b_1(G_\G) \leq {\width(\G) +1\choose 2}$.
\end{problem}

It is natural to ask  the corresponding questions for the semigroups themselves.

\begin{problem}\label{ProblemHSr1}
Let $\G$ be a numerical semigroup.
Is 
$\rho(\G) \leq \rho(\tG)$?
\end{problem}

\begin{problem}\label{ProblemHSr2}
Let $\G$ be a numerical semigroup.
Is 
$\rho(\G) \leq {\width(\G)+1 \choose 2}$?
\end{problem}

It is also natural to extend these questions  to higher Betti numbers.

\begin{problem}\label{ProblemHSb1}
Let $\G$ be a numerical semigroup.
Is 
$b_i(R_\G) \leq  b_i(R_\tG)$ for all $i$?
\end{problem}

\begin{problem}[\protect{\cite[Conjecture 1.3]{CavigliaMoscarielloSammartano}}]\label{ProblemHSb2}
Let $\G$ be a numerical semigroup, then 
$b_i(R_\G) \leq i{\width(\G)+1 \choose i+1}$
for all $i \geq 1$.
\end{problem}

Finally, one could also consider the  versions of Problems \ref{ProblemHSb1} and \ref{ProblemHSb2} for the associated graded rings.

\begin{problem}\label{ProblemHSbG1}
Let $\G$ be a numerical semigroup.
Is 
$b_i(G_\G) \leq  b_i(G_\tG)$ for all $i$?
\end{problem}

\begin{problem}\label{ProblemHSbG2}
Let $\G$ be a numerical semigroup.
Is 
$b_i(G_\G) \leq i{\width(\G)+1 \choose i+1}$ for all $i$?
\end{problem}

We now explain the relations among the  problems in this section.
Recall that $\rho(\G) = b_1(R_\G)$ and $b_i(R_\G) \leq b_i(G_\G)$.
When $\G$ is  generated by an arithmetic sequence, one has $b_i(R_\G)=b_i(G_\G)$ for all $i$,
cf.  \cite{SZN} or \cite[Proposition 2.5]{HerzogStamate},
and, moreover, 
the Betti numbers can be determined explicitly, 
cf. \cite[Theorem 4.1]{GSS}.
From this, it is possible to deduce an affirmative answer to Problem \ref{ProblemHSbG2} for semigroups generated by arithmetic sequences, cf. \cite[Proposition 2.5]{HerzogStamate}.
As a consequence, 
we have the following implications among affirmative answers to the last 8 problems:

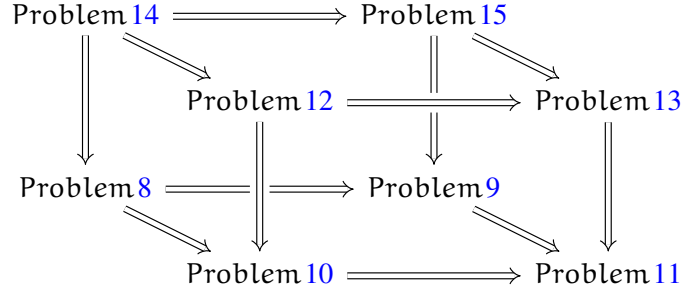
\begin{center}
\begin{tikzcd}[row sep=1.5em, column sep =.1]
Problem \, \ref{ProblemHSbG1} & & Problem \, \ref{ProblemHSbG2} &  &  \\
& Problem \, \ref{ProblemHSb1} & & Problem \, \ref{ProblemHSb2} &    \\
Problem \, \ref{ProblemHSgr1} & & Problem \, \ref{ProblemHSgr2} &  &  \\
& Problem \, \ref{ProblemHSr1} & & Problem \, \ref{ProblemHSr2} &    \\
 \arrow[Rightarrow,from=1-1,to=1-3]
 \arrow[Rightarrow,from=1-1,to=3-1]
 \arrow[Rightarrow,from=1-3,to=3-3]
 \arrow[Rightarrow,from=3-1,to=3-3]
 \arrow[Rightarrow,from=1-1,to=2-2]
 \arrow[Rightarrow,from=1-3,to=2-4]
 \arrow[Rightarrow,from=3-3,to=4-4]
 \arrow[Rightarrow,from=3-1,to=4-2]
 \arrow[Rightarrow,from=2-2,to=2-4, crossing over]
 \arrow[Rightarrow,from=2-4,to=4-4, crossing over]
 \arrow[Rightarrow,from=2-2,to=4-2, crossing over]
 \arrow[Rightarrow,from=4-2,to=4-4, crossing over]
\end{tikzcd}
\end{center}

Each of the inequalities proposed in this section, if true, would be sharp, because they are all attained by the semigroups of the form $\langle w, w+1, \ldots, 2w-1\rangle$.
Thus, any of the Problems from \ref{ProblemHSgr1} to \ref{ProblemHSbG2} gives the  following conjectural answer to Problem \ref{ProblemWidth}:
$$
\mathscr{W}(w) = {w+1 \choose 2}.
$$

Very little is known about the problems in this section.
While the results of \cite{Vu} imply that $\mathscr{W}(w)$ is finite, they give no estimate for this quantity.
In \cite{CavigliaMoscarielloSammartano}, the following  upper bounds are obtained   for the  Betti numbers of $R_\G$ and, in particular, for  $\mathscr{W}(w)$.

\begin{thm}[\protect{\cite[Theorem 1.4]{CavigliaMoscarielloSammartano}}]
Let $\G$ be a numerical semigroup. Then, 
$$b_i(R_\G) \leq {\width(\G) \choose i} (3 \mathfrak{e})^{\sqrt{2\width({\G})}}$$ 
for all $i \geq 1$, 
where $\mathfrak{e}$ denotes Euler's constant.
\end{thm}

Moreover, Problem \ref{ProblemHSb2} is settled in the first open case $e=4$ for all but finitely many values of the width.

\begin{thm}[\protect{\cite[Theorem 1.5]{CavigliaMoscarielloSammartano}}]
Let $\G$ be a numerical semigroup with $\edim(\G) = 4$.
Assume $\width(\G) \geq 40$,
then
$b_i(R_\G) \leq i{\width(\G) +1 \choose i+1}$ for all $i \geq 1$.
\end{thm}

In \cite{LTV},  the authors settle Problem \ref{ProblemHSbG1} when $\edim(\G) =3$. 
Observe that in this case, while  the problems about $R_\G$ are  solved by  Theorem \ref{ThmHerzog}, 
 the problems about $G_\G$ are nontrivial, since $G_\G$ can have arbitrarily large Betti numbers, cf. \cite[Section 3.3]{HerzogStamate}.

\begin{thm}[\protect{\cite[Theorem 1.1]{LTV}}]
Let $\G= \langle g_1, g_2, g_3 \rangle$ be a numerical semigroup with $\edim(\G) = 3$. Then, 
$b_i(G_\G) \leq b_i(G_\tG)$ for all $i$.
\end{thm}

\section{Symmetric semigroups}\label{SectionSymmetric}

A numerical semigroup $\G$ is \emph{symmetric} if for every $x \in \ZZ\setminus \G$ we have $\Frob(\G)-x \in \G$.
This is equivalent to the fact that $\type(\G)= 1$,
or that the semigroup ring $R_\G$ is Gorenstein
\cite{Kunz}.
In this section, we consider relations of symmetric semigroups.

\begin{problem}\label{ProblemSymmetric}
Let $e \in \NN$. 
Determine $\mathscr{S}(e) = \sup\big \{\rho(\G) \, \mid \, \G \text{ is symmetric and } \edim(\G) = e\big\}$.
\end{problem}

Unlike Problem \ref{ProblemWidth}, it is not known whether $\mathscr{S}(e)$ is finite for every $e\in \NN$:
in fact, this is the content of Bresinsky's problem, arguably the oldest open problem in the theory of numerical semigroups.

\begin{problem}[\protect{Bresinsky's problem \cite{Bresinsky75}}]\label{ProblemBresinsky}
Is $\mathscr{S}(e)< \infty$ for every $e \in \NN$?
\end{problem}

We have $\mathscr{S}(3) = 2$ by Theorem \ref{ThmHerzog}. 
The following theorem solves the case $e=4$.

\begin{thm}[\protect{\cite{Bresinsky75}}]\label{TheoremBresinsky4}
If $\G$ is symmetric and $\edim(\G) =4$, then $\rho(\G) =3$ or $\rho(\G) = 5$.
\end{thm}

In fact, this implies that all the Betti numbers of $R_\G$ are bounded, since we have $b_0(R_\G) = b_3(R_\G) =1$ and $b_1(R_\G) = b_2(R_\G)$ for a symmetric semigroup $\G$ with $\edim(\G) =4$,
by \eqref{EqAlternating}.

Problem \ref{ProblemBresinsky} is open for $e \geq 5$. 
The difficulty stems from the fact that it is unclear how the symmetry of $\G$ should be related to its relations.
On  one hand, the Gorenstein property  has a simple characterization in terms of  syzygies: $R_\G$ is Gorenstein if and only if its minimal free resolution is self-dual,
that is, 
$b_i(R_\G) = b_{e-1-i}(R_\G)$,
where $e = \edim(\G)$.
On the other hand, 
the self-duality of the resolution imposes no constraints on the intermediate Betti numbers, in general.
For example, it is easy to construct Gorenstein rings with Betti numbers $(1,n,n,1)$, with $n$ arbitrarily large, cf.
\cite{BE}.

In the first open case $e=5$, the following theorem gives a partial affirmative answer.

\begin{thm}[\protect{\cite{Bresinsky79}}]\label{TheoremBresinsky5}
Let $\G$ be a symmetric semigroup with $\edim(\G) =5$.
Suppose there exist four distinct indices $\{i,j,h,k\} \subseteq \{1,2,3,4,5\}$ such that $g_i+g_j = g_h+g_k$.
Then, $\rho(\G) \leq 13$.
\end{thm}

Theorem \ref{TheoremBresinsky5} is sharp, thus, we have $\mathscr{S}(5) \geq 13$.
As the author points out, the relation $g_i+g_j = g_h+g_k$
is interesting because it appears in many families of numerical semigroups with unbounded number of relations.

\section{Almost symmetric semigroups}\label{SectionAlmostSymmetric}

A numerical semigroup $\G$ is \emph{almost symmetric} if,
for every $x \in \ZZ\setminus \G$,  either  $\Frob(\G)-x \in \G$ or both $x$ and $ \Frob(\G)-x$ are pseudo-Frobenius numbers of $\G$.
Almost symmetric semigroups of type 1 are precisely the symmetric semigroups.

\begin{problem}\label{ProblemAlmostSymmetric}
Let $e \in \NN$. 
Determine $\mathscr{A}(e) = \sup\big \{\type(\G) \, \mid \, \G \text{ is almost symmetric and } \edim(\G) = e\big\}$.
\end{problem}

Again, we have $\mathscr{A}(3) = 2$ by Theorem \ref{ThmHerzog}, and
the following theorem solves the case $e=4$.

\begin{thm}[\protect{\cite[Theorem 1]{Moscariello16}}]\label{TheoremMoscariello4}
We have $\mathscr{A}(4) = 3$. 
\end{thm}

As in Section \ref{SectionSymmetric},  it is not known whether $\mathscr{A}(e)$ is finite for every $e\in \NN$.

\begin{problem}[\protect{ \cite[Question 13]{Moscariello16}}]\label{ProblemMoscariello}
Is $\mathscr{A}(e)< \infty$ for every $e \in \NN$?
\end{problem}

The answer is affirmative when $e = 5$.

\begin{thm}[\protect{\cite[Theorem 2]{Moscariello23}}]
Let $\G$ be an almost symmetric  semigroup with $\edim(\G) = 5$.
Then, $\type(\G) \leq 473$.
\end{thm}

The exact value of 
$\mathscr{A}(5) $ is unknown,
though it is expected to be much smaller than 473.

Next, we discuss two generalizations of almost symmetric semigroups. 
A numerical semigroup $\G$ is \emph{nearly Gorenstein} if the following condition holds:
for every generator $g_i$ of $\G$ there exists a pseudo-Frobenius number $p$ such that $g_i+p-q \in \G$ for every pseudo-Frobenius number $q$ \cite{MoscarielloStrazzanti}.
Every almost symmetric semigroup is nearly Gorenstein, but the converse is not true.
Theorem \ref{TheoremMoscariello4} was later extended to nearly Gorenstein semigroups:

\begin{thm}[\protect{\cite[Theorem 2.4]{MoscarielloStrazzanti}}]
Let $\G$ be a nearly Gorenstein   semigroup with $\edim(\G) = 4$.
Then, $\type(\G) \leq 3$.
\end{thm}

\begin{problem}[\protect{\cite[Question 4.2]{MoscarielloStrazzanti}}]\label{ProblemMoscarielloStrazzanti}
Let $\G$ be a nearly Gorenstein  semigroup with $\edim(\G) =5$.
Is it true that $\type(\G) \leq 5$,  and that the equality is attained only if $\G$ is almost symmetric?
\end{problem}

In particular, this would imply that $\mathscr{A}(5) =5$.
We observe that $\type(\G) \leq \edim(\G)$ may fail  for almost symmetric semigroups with $\edim(\G) \geq 6$. 
In fact, the difference $\type(\G) -\edim(\G)$ can be arbitrarily large, cf. \cite[p. 277]{Moscariello16}.

A numerical semigroup $\G$ has a \emph{canonical reduction} if $g_1+\Frob(\G)-q \in \G$ for  
every pseudo-Frobenius number $q$ \cite{MoscarielloStrazzanti}.
Every nearly Gorenstein  semigroup has a canonical reduction, but the converse is not true.

\begin{problem}[\protect{\cite[Question 4.4]{MoscarielloStrazzanti}}]
Let $\G$ be a numerical semigroup with canonical reduction and $\edim(\G) =4$.
Is $\type(\G) \leq 4$?
\end{problem}

One can generalize these problems to all the Betti numbers.

\begin{problem}[\protect{\cite[Question 3.7]{Stamate}}]\label{ProblemStamate}
Let $\G$ be an almost symmetric or nearly Gorenstein numerical semigroup.
Are there upper bounds for the Betti numbers $b_i(R_\G)$ in terms of  $\edim(\G)$?
\end{problem}

The case $i=1$ of Problem \ref{ProblemStamate} is of particular interest, as it generalizes  Bresinsky's Problem \ref{ProblemBresinsky}.
Theorem \ref{TheoremBresinsky4} was extended to almost symmetric semigroups in \cite{Eto} and, independently, in \cite{HerzogWatanabe}.

\begin{thm}[\protect{\cite{Eto,HerzogWatanabe}}]\label{TheoremEtoHW}
If $\G$ is an almost symmetric  semigroup with $\edim(\G) = 4$, 
then $\rho(\G) \leq 7$.
\end{thm}

As in Section \ref{SectionSymmetric}, by \eqref{EqAlternating}, this implies that all the Betti numbers are bounded.

The analogue of Theorem \ref{TheoremEtoHW}  for nearly Gorenstein semigroups is not known.

\begin{problem}[\protect{\cite[Question 4.3]{MoscarielloStrazzanti}}]
Let $\G$ be a nearly Gorenstein numerical semigroup  which is not almost symmetric. 
Assume that  $\edim(\G)=4$.
Is it true that $\rho(\G) \in \{5,6\}$?
\end{problem}

\section{Uniform boundedness}

Inspired by the questions and results of Sections \ref{SectionEdimMult}, \ref{SectionSymmetric} and \ref{SectionAlmostSymmetric},
we propose a general problem on the boundedness of Betti numbers.

\begin{problem}\label{ProblemBoundedness}
Let  $\mathcal{F}$ be a family of numerical semigroups with fixed embedding dimension $e$.
Define $\mathscr{B}_i(\mathcal{F}) = \sup \big\{ b_i(R_\G) \, \mid \, \G \in \mathcal{F}\big\}$ for $i = 1,\ldots,e-1$.
Is it true that the $\mathscr{B}_i(\mathcal{F}) $ are either all finite or all infinite?
\end{problem}

All families of semigroups for which all the Betti numbers are known follow this behavior, see for example
\cite{Arslan,Bresinsky75,Bresinsky75b,Eto,FGH,Gluing,GSS,HerzogWatanabe,MSS,Stamate}.

By Theorem \ref{ThmHerzog}, the first interesting case of Problem \ref{ProblemBoundedness} is $e=4$,
and the problem is already open in this case.
Recall that, by \eqref{EqAlternating}, we have $b_2(R_\G) = \rho(\G) + \type(\G) -1$ when $\edim(\G) =4$,
thus, it suffices to show that $\rho(\G)$ is bounded if and only if $\type(\G)$ is bounded.
One direction is settled in the following theorem.

\begin{thm}[\protect{\cite[Theorem 7]{Bresinsky88}}]
Let $\G$ be a  semigroup with $\edim(\G) = 4$.
Then, 
$\rho(\G) \leq 4+9\type(\G)$.
\end{thm}

It would be therefore interesting, and perhaps feasible, to investigate the other direction, which would settle Problem \ref{ProblemBoundedness} in the case of embedding dimension 4.

\begin{problem}
Let $\G$ be a  semigroup with $\edim(\G) = 4$.
Is $\type(\G)$ bounded by a function of 
$\rho(\G)$?
\end{problem}

We conclude 
with the following compelling question, which relates the first syzygies and the last syzygies of $R_\G$, 
thus fitting the theme of this section.
Let $p$ be a pseudo-Frobenius number of a  semigroup $\G = \langle g_1, \ldots, g_e\rangle$.
An \emph{RF-matrix} (row-factorization matrix) of $p$ is an $e\times e$ matrix $(a_{ij})$ such that $a_{ii}=-1$ for all $i$, $a_{ij} \in \NN$ whenever $i \ne j$, and $\sum_{j=1}^e a_{ij}g_j = p$ for all $i$ \cite{Moscariello16}. 
Given two rows $\mathbf{a}_i,\mathbf{a}_j \in \ZZ^e$ of an RF-matrix $(a_{ij})$, 
write $\mathbf{a}_i-\mathbf{a}_j = \mathbf{a}^+-\mathbf{a}^-$ with $\mathbf{a}^+,\mathbf{a}^-\in \NN^e$.
In the notation of Section \ref{SectionSetup}, we have $(\mathbf{a}^+,\mathbf{a}^-) \in \ker(\varphi)$,
and $(\mathbf{a}^+,\mathbf{a}^-)$ is called an RF-relation of $\G$.

\begin{solved}[\protect{\cite[Question 3.6]{HerzogWatanabe}}] \label{ProblemRFHW}
Are all minimal relations of a numerical semigroups RF-relations?
\end{solved}

See \cite{BGS} for some partial results in small embedding dimension. 

Problem \ref{ProblemRFHW}
 was solved in \cite[Corollary 3.5]{EMMNSW}.

\section{Complete intersections}\label{SectionCompleteIntersection}

So far, we have mostly focused on upper bounds for the relations of a numerical semigroup. 
The picture for lower bounds is much clearer: we always have $\rho(\G) \geq \edim(\G)-1$,
and this inequality is sharp.
A numerical semigroup attaining  equality is called a \emph{complete intersection}.
Complete intersection semigroups are always symmetric, but the converse is not true if $\edim(\G) \geq 4$.
A semigroup $\G$ is a complete intersection if and only if the semigroup ring $A_\G$ (equivalently, $R_\G$) is a complete intersection, and, in this case,  $A_\G$ is resolved by a Koszul complex
and $b_i(A_\G) = {\edim(\G) -1 \choose i}$.

Complete intersections are, 
in some sense,
the best-behaved numerical semigroups,
and they have been the subject of extensive research \cite{AGS, DMS,GMGSOT}.
In particular, they can be  characterized recursively in terms of gluing \cite{Delorme}.
Nevertheless,
we still do not fully understand  the implications of the complete intersection condition,  and
there are  important unsolved problems on this topic.

The Hilbert series of $A_\G$ is the power series $H_\G(z) = \sum_{\gamma \in \G} z^\gamma \in \ZZ\llbracket z\rrbracket$.
It is a rational function,
in fact,
one has $H_\G(z) = \frac{P_\G(z)}{1-z} $ where 
$P_\G(z) = 1 +(z-1) \sum_{\gamma \in \NN \setminus \G} z^\gamma$.
The polynomial $P_\G(z)$ is monic, with integer coefficients, and degree equal to $\Frob(\G)+1$.
A semigroup $\G$ is called \emph{cyclotomic} if $P_\G(z)$ is a product of cyclotomic polynomials. 
Cyclotomic semigroups are always symmetric, but the converse is not true \cite{CGSM,SS}.

\begin{problem}[\protect{\cite[Conjecture 1]{CGSM}}]\label{ProblemCyclotomic}
A numerical semigroup is a complete intersection if and only if it is cyclotomic.
\end{problem}

The forward direction of this conjecture is easy, 
since  the Hilbert series of a positively 
$\NN$-graded complete intersection is a product of factors of the form $(1-z^a)^b$ with $a,b\in\ZZ$.
But the backward direction is striking: 
in general, 
the Hilbert series of a graded algebra is a much coarser invariant than its free resolution,
and it is not possible to detect the complete intersection property from the Hilbert series alone.
For instance, there exists 
a normal Gorenstein toric domain with the same Hilbert series as a complete intersection, but which is not a complete intersection \cite[Example 3.9]{Stanley}.
However, there are positive results in this spirit,
for instance,
notably, Stanley's theorem 
characterizing the Gorenstein property for 
positively $\NN$-graded
 Cohen-Macaulay domains in terms of a symmetry in the Hilbert series
\cite[Theorem 4.4]{Stanley}.
See \cite{BM} for more discussions on this theme.

Problem \ref{ProblemCyclotomic} remains widely open, 
despite some effort \cite{BHPM,CGSHPM,SS}.
It is open already when $\edim(\G) = 4$.

The Hilbert function of $G_\G=\gr(R_\G)$ is the numerical function $\mathrm{HF}_{G_\G}:\NN\to \NN$ defined as
$$
\mathrm{HF}_{G_\G}(j) = \dim_\kk [G_\G]_j= \dim_\kk \big({\mm^{j}}/{\mm^{j+1}}\big).
$$

\begin{problem}[\protect{Rossi's Problem \cite[Problem 4.9]{Rossi}}]\label{ProblemRossi}
Let $\G$ be a complete intersection numerical semigroup.
Is the Hilbert function  $\mathrm{HF}_{G_\G}$ non-decreasing?
\end{problem}

The answer is obviously  affirmative when  $G_\G$ is  Cohen-Macaulay.
However, the fact that the local ring $R_\G$ is a complete intersection does not imply that  its 
associated graded ring  $G_\G$ is a complete intersection, or even Cohen-Macaulay.
This is what makes Problem \ref{ProblemRossi}  interesting, and also challenging: 
it is unclear how the assumption on the relations of $R_\G$ should constrain the Hilbert function of $G_\G$.

When $\edim(\G) \leq 3$, the answer is affirmative, even without any  assumption on $\G$ \cite[Theorem 4.6]{Rossi}.
Problem \ref{ProblemRossi} was originally stated for symmetric semigroups (in fact, for Gorenstein one-dimensional local rings), but the answer is negative in this generality \cite{OST}.
But the most interesting case is that of complete intersections, which remains open.
Some partial results were obtained in \cite{AM,ASS}.

We conclude this section with another  question about the structure of complete intersection numerical semigroups.

\begin{solved}[\protect{\cite[Question 7.5]{EHOOPK}}]\label{ProblemCIem}
Characterize the pairs $(e,m)$ for which there exists a complete intersection numerical semigroup with $\edim(\G) = e$ and $\mult(\G) = m$.
\end{solved}

In \cite[Proposition 5]{AGS}, it is proved that complete intersection semigroups satisfy $\mult(\G) \geq 2^{\edim(\G)-1}$,
and that this bound is attained.

Problem \ref{ProblemCIem} was actually solved long before this problem list appeared, 
in \cite[Theorem 1]{Watanabe}.

The version of Problem \ref{ProblemCIem} for symmetric semigroups was solved in \cite[Theorem 5]{Rosales01},
where it is proved that there exists a symmetric numerical semigroup
with $\edim(\G) = e$ and $\mult(\G) = m$ if and only if $2 \leq e \leq m-1$ or $(e,m) = (1,1), (2,2)$.

We refer the reader to \cite[Section 6]{AGS}, 
\cite[Section 4]{GSL},
and 
\cite[Sections 3 and 5]{GMGSOT} for more open questions on complete intersection numerical semigroups.

\section{Frobenius number}

When considering the graded semigroup ring $A_\G$, 
instead of the local semigroup ring $R_\G$,
one can study invariants related to the graded syzygies of $A_\G$ and, more generally, to the graded components of (co)homological functors associated to $A_\G$.
Some  care is required, since $A_\G$ is a positively $\NN$-graded $\kk$-algebra, but it is not standard graded;
see \cite{DS} for a treatment at this level of generality.
The most important  invariant of graded free resolutions is perhaps the \emph{Castelnuovo-Mumford regularity}.
It  can be computed in two ways, cf. \cite[Section 3]{DS}:
\begin{equation}\label{EqRegularity}
\reg(A_\G) 
= \max_{0 \leq i \leq e} \big\{ j +i \, \mid \, [H^i_\mathfrak{m}(A_\G)]_j \ne 0\big\}
= \max_{0 \leq i \leq e} \big\{ j -i \, \mid \, b_{i,j}(A_\G) \ne 0\big\} - \sum_{h=1}^e g_h + e.
\end{equation}
Here,  $\mathfrak{m}$ denotes the homogeneous maximal ideal of $A_\G$,
and $H^i_\mathfrak{m}(A_\G)$ denotes the $i$th local cohomology module of $A_\G$.
Using the first equation in \eqref{EqRegularity} and the fact that $A_\G$ is Cohen-Macaulay of dimension one, it is not hard to see that 
\begin{equation}\label{EqRegularityFrobenius}
\reg(A_\G) = \Frob(\G)+1.
\end{equation}
This allows to connect the vast literature on the Frobenius problem \cite{RA} to the topic of free resolutions.

One of the most famous open problems about the Frobenius number is Wilf's Conjecture. 
In order to state it, we define the integer $\eta(\G) = \big|\{ \gamma \in \G \, \mid \, \gamma < \Frob(\G)\}\big|$.

\begin{problem}[Wilf's Conjecture \protect{\cite{Wilf}}]\label{ProblemWilf}
Let $\G$ be a numerical semigroup,
then $\Frob(\G) < \edim(\G) \eta(\G)$.
\end{problem}

Wilf's conjecture has been tackled in many papers, and partial progress 
in several directions has been achieved.
For example, we mention \cite{Eliahou}, where tools from commutative algebra are used to settle the conjecture under the assumption that $\Frob(\G) < 3 \,\mult(\G)$.
For an extensive overview of the state of the art, we refer to \cite{Delgado}.
Here, we limit ourselves to formulating 
a weaker version of Wilf's conjecture.

By \eqref{EqRegularity} and \eqref{EqRegularityFrobenius}, Wilf's inequality can be regarded as an upper bound on the degrees of the syzygies of the semigroups ring.
Restricting to the first syzygies, that is, to the defining ideal of $A_\G$, we obtain the following natural question, which may be more tractable than Problem \ref{ProblemWilf}.

\begin{problem}[Weak Wilf's Conjecture]\label{ProblemWilfIdeal}
Let $f$ be a minimal  generator of the toric ideal of $A_\G$.
Then,
$$
\deg(f) 
\leq \edim(\G) (\eta(\G) -1) + \sum_{i=1}^{\edim(\G)} g_i +1. 
$$
\end{problem}

\addtocontents{toc}{\protect\setcounter{tocdepth}{0}}

\section*{Acknowledgments}
We would like to thank 
Pedro A. García-Sánchez,
Vincenzo Micale,
Dumitru Stamate,
Francesco Strazzanti,
and
Thanh Vu
for comments on a preliminary version of this paper.

\section*{Funding}
Moscariello was supported by 
the grant “Proprietà locali e globali di
anelli e di varietà algebriche” PIACERI 2020–22, Università degli Studi di Catania.
Sammartano was  supported by the grant PRIN 2020355B8Y
{\em Square-free Gr\"obner degenerations, special varieties and related topics}
and by the INdAM – GNSAGA Project CUP E55F22000270001.

\addtocontents{toc}{\protect\setcounter{tocdepth}{1}}


\begin{thebibliography}{aaaaaaaaaa}

\bibitem[A99]{Arslan}
F. Arslan, 
\href{https://doi.org/10.1090/S0002-9939-99-05229-6}
{\emph{Cohen-Macaulayness of tangent cones}},
 Proc. Am. Math. Soc. {\bf 128} (2000), 2243--2251.

\bibitem[AM07]{AM}
F. Arslan, P. Mete,
\href{https://doi.org/10.1090/S0002-9939-07-08793-X}
{\emph{Hilbert functions of Gorenstein monomial curves}}, 
Proc. Amer. Math. Soc. {\bf 135}
(2007), 1993--2002.

\bibitem[ASS13]{ASS}
F. Arslan, N. Sipahi, N. Şahin, 
\href{https://doi.org/10.1016/j.jsc.2013.03.002}
{\emph{Monomial curve families supporting Rossi’s conjecture}}, 
J. Symbolic
Comput. {\bf 55} (2013), 10--18.

\bibitem[AGS13]{AGS}
A. Assi, P. A. García-Sánchez,
\href{https://doi.org/10.1007/s00200-013-0186-z}
{\emph{Constructing the set of complete intersection numerical semigroups with a given Frobenius number}},
 Appl. Algebra  Eng. Commun.  Comput. {\bf 24} (2013), 133--148.


\bibitem[B06]{Barucci}
V. Barucci, 
\href{https://link.springer.com/content/pdf/10.1007/978-0-387-36717-0_3.pdf}
{\emph{Numerical semigroup algebras}}.
 Multiplicative Ideal Theory in Commutative Algebra: A Tribute to the Work of Robert Gilmer. Boston, MA: Springer US (2006), 39--53.

\bibitem[BGS22]{BGS}
O. P. Bhardwaj, K. Goel,  I. Sengupta,
\href{https://doi.org/10.1142/S0218196722500552}
{\emph{On row-factorization relations of certain numerical semigroups}},
Int. J. Algebra Comput. {\bf  32} (2022), 1275--1305.


\bibitem[BM04]{BM}
R. B\o gvad, T. Meyer,
\href{https://doi.org/10.1016/j.jsc.2004.06.001}
{\emph{On algorithmically checking whether a Hilbert series comes from a complete intersection}},
J. Symb. Comput. {\bf 38} (2004), 1487--1506.

\bibitem[BHPM21]{BHPM}
A. Borzì, A. Herrera-Poyatos, P. Moree,
\href{https://doi.org/10.1007/s00233-021-10197-8}
{\emph{Cyclotomic numerical semigroup polynomials with at most two irreducible factors}},
Semigroup Forum. {\bf 103} (2021), 812--828.

\bibitem[BE77]{BE}
D. A. Buchsbaum, D. Eisenbud,
\href{https://doi.org/10.2307/2373926}
{\emph{Algebra structures for finite free resolutions, and some structure theorems for ideals of codimension 3}},
Amer. J. Math. {\bf 99} (1977), 447--485.

\bibitem[B75]{Bresinsky75}
H. Bresinsky, 
\href{https://doi.org/10.1007/BF01170309}
{\emph{Symmetric semigroups of integers generated by 4 elements}},
Manuscripta Math.  {\bf 17} (1975), 205--219.


\bibitem[B75b]{Bresinsky75b}
H. Bresinsky, 
\href{https://doi.org/10.1090/S0002-9939-1975-0389912-0}
{\emph{On prime ideals with generic zero $x_i = t^{n_i}$}},
Proc. Amer. Math. Soc. {\bf 47} (1975), 
329--332.


\bibitem[B79]{Bresinsky79}
H. Bresinsky, 
\href{https://doi.org/10.1007/BF01303625}
{\emph{Monomial gorenstein ideals}},
Manuscripta Math.  {\bf 29} (1979), 159--181.
159-181.

\bibitem[B88]{Bresinsky88}
H. Bresinsky, 
\href{https://doi.org/10.1007/BF01170309}
{\emph{Binomial generating sets for monomial curves, with applications in $\mathbb{A}^4$}},
 Rend. Sem. Mat. Univ. Politec. Torino {\bf 46} (1988), 353--370.

\bibitem[BGT02]{BGT}
W. Bruns, J. Gubeladze, N. V. Trung, 
\href{https://doi.org/10.1007/s002330010099}
{\emph{Problems and algorithms for affine semigroups}},
Semigroup Forum {\bf 64} (2002), 180--212.

\bibitem[BH97]{BrunsHerzog}
W. Bruns, J. Herzog,
\href{https://doi.org/10.1016/S0022-4049(97)00051-0}
{\emph{Semigroup rings and simplicial complexes}},
J. Pure Appl. Algebra
{\bf 122} (1997), 185--208.

\bibitem[CM91]{CampilloMarijuan}
A. Campillo, C. Marijuan,
\href{https://jtnb.centre-mersenne.org/item/JTNB_1991__3_2_249_0/}
{\emph{Higher order relations for a numerical semigroup}},
 Sém. Théor. Nombres Bordeaux {\bf 3} (1991), 249--260.

\bibitem[CMS24]{CavigliaMoscarielloSammartano}
G. Caviglia, A. Moscariello, A. Sammartano,
\href{https://doi.org/10.1090/proc/16862}
{\emph{Bounds for syzygies of monomial curves}},
 Proc. Amer. Math. Soc.
{\bf 152} (2024), 3665--3678.


\bibitem[CGSM16]{CGSM}
E. A. Ciolan, P. A. García-Sánchez, P. Moree,
\href{https://doi.org/10.1137/140989479}
{\emph{Cyclotomic numerical semigroups}},
 SIAM J. Discrete Math. {\bf 30} (2016), 650--668.

\bibitem[CGSHPM22]{CGSHPM}
E. A. Ciolan, P. A. García-Sánchez, A. Herrera-Poyatos, P. Moree,
\href{https://doi.org/10.1016/j.disc.2022.112820}
{\emph{Cyclotomic exponent sequences of numerical semigroups}},
Discrete Math.{\bf 345} (2022), 112820.

\bibitem[DS08]{DS}
G. Dalzotto, E. Sbarra,
\href{https://arxiv.org/abs/math/0506333}
{\emph{On non-standard graded algebras}},
Toyama Math. J. {\bf 31} (2008), 33--57.


\bibitem[DMS14]{DMS}
M. D’Anna, V.  Micale, A. Sammartano,
\href{https://doi.org/10.1007/s00233-013-9547-y}
{\emph{Classes of complete intersection numerical semigroups}},
 Semigroup Forum {\bf 88} (2014), 453--467.

\bibitem[D20]{Delgado}
M. Delgado, 
\href{https://doi.org/10.1007/978-3-030-40822-0_4}
{\emph{Conjecture of Wilf: a survey}},
 Numerical Semigroups: IMNS 2018 (2020), 39--62.

\bibitem[D76]{Delorme}
C. Delorme,
\href{http://www.numdam.org/item/10.24033/asens.1307.pdf}
{\emph{Sous-monoïdes d’intersection complète de $ N$}},
Ann. Sci. Éc. Norm. Supér. {\bf  9} (1976), 
145--154.

\bibitem[E95]{Eisenbud}
D.  Eisenbud,
\href{https://www.springer.com/gp/book/9780387942681}
{\emph{Commutative Algebra: with a view toward algebraic geometry}}, 
Vol. 150. Springer Science \& Business Media (1995).

\bibitem[EGV91]{EGV}
J. Elias, A. V. Geramita, G. Valla,
\href{https://doi.org/10.1090/conm/159}
{\emph{On the Cohen-Macaulay type of perfect ideals}},
Contemp. Math. {\bf 159} (1994), 41--49.

\bibitem[ERV91]{ERV}
J. Elias, L. Robbiano, G. Valla,
\href{https://doi.org/10.1017/S0027763000003640}
{\emph{Number of generators of ideals}},
Nagoya Math. J. {\bf 123} (1991), 39--76.

\bibitem[E18]{Eliahou}
S. Eliahou, 
\href{https://ems.press/content/serial-article-files/32312}
{\emph{Wilf's conjecture and Macaulay's theorem}},
J. Eur. Math. Soc. {\bf 20} (2018), 2105--2129.

\bibitem[EHOOPK24]{EHOOPK}
C. Elmacioglu, K. Hilmer, C. O'Neill, M. Okandan, H. Park-Kaufmann,
\href{https://doi.org/10.5802/alco.354}
{\emph{On the cardinality of minimal presentations of numerical semigroups}},
Algebr. Comb. {\bf 7} (2024), 753--771.


\bibitem[E17]{Eto}
K. Eto, 
\href{https://doi.org/10.1016/j.jalgebra.2017.05.044}
{\emph{Almost Gorenstein monomial curves in affine four space}},
J. Algebra {\bf 488} (2017), 362--387.

\bibitem[EMMNSW26]{EMMNSW}
K. Eto, N.  Matsuoka, A. Moscariello, T. Numata, A. Sammartano, K. Watanabe, 
\href{https://arxiv.org/abs/2609.24880}
{\emph{The structure of almost symmetric almost complete intersection numerical semigroups}},
 arXiv 2609.24880 (2026).

\bibitem[FGH86]{FGH}
R. Fröberg, C. Gottlieb, R. H\"aggkvist. 
\href{https://doi.org/10.1007/BF02573091}
{\emph{On numerical semigroups}}
Semigroup Forum {\bf 35} (1986), 63--83.


\bibitem[GMGSOT24]{GMGSOT}
I. García-Marco, P. A. García-Sánchez, I. Ojeda, C. Tatakis, 
\href{https://doi.org/10.1016/j.jpaa.2023.107551}
{\emph{Universally free numerical semigroups}},
J. Pure Appl. Algebra {\bf 228} (2024), 107551.

\bibitem[GSL13]{GSL}
P. A. Garcia-Sanchez, M. J. Leamer,
\href{https://doi.org/10.1016/j.jalgebra.2013.06.007}
{\emph{Huneke–Wiegand conjecture for complete intersection numerical semigroup rings}}, J. Algebra {\bf 391} (2013), 114--124.

\bibitem[GSS13]{GSS}
P. Gimenez, I. Sengupta, I. Srinivasan, 
\href{https://doi.org/10.1016/j.jalgebra.2013.04.026}
{\emph{Minimal graded free resolutions for monomial curves defined by arithmetic sequences}}, 
J. Algebra {\bf 338} (2013), 294--310.

\bibitem[GS19]{Gluing}
P. Gimenez, H. Srinivasan. 
\href{https://doi.org/10.1016/j.jpaa.2018.06.010}
{\emph{The structure of the minimal free resolution of semigroup rings obtained by gluing}},
J. Pure Appl. Algebra  {\bf 223} (2019), 1411--1426.

\bibitem[GS20]{GS20}
P. Gimenez, H. Srinivasan. 
\href{https://doi.org/10.1007/978-3-030-40822-0_8}
{\emph{Syzygies of Numerical Semigroup Rings, a Survey Through Examples}}.
In: Numerical Semigroups. Springer INdAM Series, vol 40. Springer, Cham. 
(2020).

\bibitem[M2]{M2}
D.R. Grayson, M.E. Stillman, 
\href{http://www.math.uiuc.edu/Macaulay2/}{\textit{Macaulay 2, a software system for research in algebraic geometry}}.

\bibitem[H70]{Herzog}
J. Herzog, 
\href{https://doi.org/10.1007/BF01273309}
{\emph{Generators and relations of abelian semigroups and semigroup rings}},
Manuscr. Math.
{\bf 3} (1970), 175--193.

 \bibitem[HH11]{HerzogHibi} 
J. Herzog, T. Hibi, 
 \href{https://doi.org/10.1007/978-0-85729-106-6}
 {\emph{Monomial Ideals}}, 
 Springer, London (2011).


\bibitem[HS14]{HerzogStamate}
J. Herzog, D. I. Stamate, 
\href{https://doi.org/10.1016/j.jalgebra.2014.07.008}
{\emph{On the defining equations of the tangent cone of a numerical semigroup ring}},
J. Algebra {\bf 418} (2014), 8--28.

\bibitem[HW19]{HerzogWatanabe}
J. Herzog, K. Watanabe, 
\href{https://doi.org/10.1007/s00233-019-10007-2}
{\emph{Almost symmetric numerical semigroups}},
 Semigroup Forum {\bf 98} (2019),
 589--630.


\bibitem[K70]{Kunz}
E. Kunz,
\href{https://community.ams.org/journals/proc/1970-025-04/S0002-9939-1970-0265353-7/S0002-9939-1970-0265353-7.pdf}
{\emph{The value-semigroup of a one-dimensional Gorenstein ring}},
Proc. Amer. Math. Soc. {\bf 25} (1970), 748--751.

\bibitem[LTV24]{LTV}
N. P. H. Lan, N. C. Tu, T. Vu,
\href{https://doi.org/10.1007/s40306-024-00546-4}
{\emph{Betti numbers of the tangent cones of monomial space curves}},
Acta Math. Viet. {\bf 49} (2024), 347--365.

\bibitem[MSS19]{MSS}
R. Mehta, J. Saha, I. Sengupta. 
\href{https://doi.org/10.1142/S0219498819501433}
{\emph{Betti numbers of Bresinsky’s curves in $\mathbb{A}^4$}},
J. Algebra Appl. {\bf 18} (2019), 1950143.

\bibitem[M16]{Moscariello16}
A. Moscariello,
\href{https://doi.org/10.1016/j.jalgebra.2016.02.019}
{\emph{On the type of an almost Gorenstein monomial curve}},
J. Algebra {\bf 456} (2016), 266--277.


\bibitem[M23]{Moscariello23}
A. Moscariello,
\href{https://doi.org/10.1080/00927872.2022.2131803}
{\emph{On the boundedness of the type of an almost Gorenstein monomial curve in $\mathbb{A}^5$}},
Commun. Algebra {\bf 51} (2023), 1179--1185.

\bibitem[MS25]{MoscarielloSammartano}
A. Moscariello, A. Sammartano,
\href{https://arxiv.org/abs/2405.19810}
{\emph{On minimal presentations of numerical monoids}},
 Bull. Lond. Math. Soc. {\bf 57}  (2025), 878--894.
 

\bibitem[MS21]{MoscarielloStrazzanti}
A. Moscariello, F. Strazzanti,
\href{https://doi.org/10.1007/s00009-021-01761-1}
{\emph{Nearly Gorenstein vs almost Gorenstein affine monomial curves}}
Mediterr. J. Math. {\bf 18} (2021), 127.

\bibitem[OST17]{OST}
A. Oneto, F. Strazzanti, G. Tamone,
\href{https://doi.org/10.1016/j.jalgebra.2017.05.038}
{\emph{One-dimensional Gorenstein local rings with decreasing Hilbert function}},
J. Algebra {\bf 489} (2017), 91-114.


\bibitem[RA05]{RA}
J. L. Ramírez-Alfonsín,
\href{https://doi.org/10.1093/acprof:oso/9780198568209.001.0001}
{\emph{The diophantine Frobenius problem}},
OUP Oxford (2005).

\bibitem[R01]{Rosales01}
J. Rosales, 
\href{https://doi.org/10.1090/S0002-9939-01-05819-1}
{\emph{Symmetric numerical semigroups with arbitrary multiplicity and embedding dimension}},
Proc. Amer. Math. Soc. {\bf 129} (2001), 2197--2203.


\bibitem[RGS98]{RosalesGarciaSanchez98}
J.C. Rosales , P. A. García-Sánchez,
\href{https://doi.org/10.1006/jabr.1997.7341}
{\emph{On numerical semigroups with high embedding dimension}},
J. Algebra {\bf 203} (1998), 567--578.


\bibitem[RGS09]{RosalesGarciaSanchezBook}
J.C. Rosales , P. A. García-Sánchez,
\href{https://doi.org/10.1007/978-1-4419-0160-6}
{\emph{Numerical semigroups}},
 Vol. 20. New York: Springer (2009).

\bibitem[R11]{Rossi}
M.E. Rossi, 
\href{https://doi.org/10.1090/conm/555}
{\emph{Hilbert functions of Cohen–Macaulay local rings}}, 
 Contemp. Math. {\bf 555} (2011), 173--200.


\bibitem[SS18]{SS}
M. Sawhney, D. Stoner,
\href{https://doi.org/10.1137/17M1138479}
{\emph{On symmetric but not cyclotomic numerical semigroups}},
SIAM J. Discrete Math. {\bf 32} (2018), 1296--1304.

\bibitem[SZN12]{SZN}
L. Sharifan, R. Zaare-Nahandi,
\emph{A class of monomial curves of homogeneous type},
 Extended Abstracts of the 22nd Iranian Algebra Seminar, Hakim Sabzevari University,
Sabzevar, Iran (2012) pp. 255–258.

\bibitem[S18]{Stamate}
D. Stamate,
\href{https://doi.org/10.1007/978-3-319-90493-1_8}
{\emph{Betti numbers for numerical semigroup rings}},
 Multigraded Algebra and Applications (2016), 133--157.


\bibitem[S78]{Stanley}
R. P. Stanley,
\href{https://doi.org/10.1016/0001-8708(78)90045-2}
{\emph{Hilbert functions of graded algebras}},
Adv. Math. {\bf 28} (1978), 57--83.

\bibitem[V94]{Valla}
G. Valla,
\href{http://www.numdam.org/item/CM_1994__91_3_305_0.pdf}
{\emph{On the Betti numbers of perfect ideals}},
Comp. Math. {\bf 91} (1994), 305--319.

\bibitem[V14]{Vu}
T. Vu, 
\href{https://doi.org/10.1016/j.jalgebra.2014.07.007}
{\emph{Periodicity of Betti numbers of monomial curves}},
J. Algebra {\bf 418} (2014), 66--90.

\bibitem[W73]{Watanabe}
K. Watanabe,
\href{https://doi.org/10.1017/S0027763000015312}
{\emph{Some examples of one dimensional Gorenstein domains}},
Nagoya Math. J. {\bf 49}  (1973), 101--109.

\bibitem[W78]{Wilf}
H. S. Wilf, 
\href{https://doi.org/10.1080/00029890.1978.11994639}
{\emph{A circle-of-lights algorithm for the “money-changing problem”}},
Amer. Math. Monthly {\bf 85} (1978), 562--565.

\end{thebibliography}
\end{document}